\documentclass[12pt]{article}
\usepackage{nicefrac}
\usepackage{hyperref}
\usepackage{amssymb,mathtools,amsthm, amsmath, mathabx}
\usepackage[all,cmtip]{xy}
\usepackage{xfrac,mathdots}
\usepackage{rotating, setspace}
\usepackage{fullpage,comment,authblk,stackrel}
\usepackage[small, labelfont = bf]{caption}
\usepackage{tikz-cd}
\usepackage[numbers,sort]{natbib}

\newcommand{\blocktheorem}[1]{%
  \csletcs{old#1}{#1}% Store \begin
  \csletcs{endold#1}{end#1}% Store \end
  \RenewDocumentEnvironment{#1}{o}
    {\par\addvspace{1.5ex}
     \noindent\begin{minipage}{\textwidth}
     \IfNoValueTF{##1}
       {\csuse{old#1}}
       {\csuse{old#1}[##1]}}
    {\csuse{endold#1}
     \end{minipage}
     \par\addvspace{1.5ex}}
}

\hypersetup{
  colorlinks   = true, %Colours links instead of ugly boxes
  urlcolor     = blue, %Colour for external hyperlinks
  linkcolor    = blue, %Colour of internal links
  citecolor   = blue %Colour of citations
}

\newcommand{\BD}       {{{ZB}}}
\newcommand{\com}       {{\text{c}}\,}
\newcommand{\Int}       {{\text{Int}}\,}
\newcommand{\Vect}      {{\mathsf{Vec}}}
\newcommand{\Gal}       {{\text{Gal}}}
\newcommand{\Z}        	{{\mathbb Z}}
\newcommand{\st}        {{\text{st}\,}}
\newcommand{\Hom}       {{\text{Hom}}}
\newcommand{\field}		{\mathbb{F}}
\newcommand\define[1]	{{\bf{#1}}}	
\newcommand{\image}		{\text{im}\,}

\newtheoremstyle{amit}% name
{7pt}% Space above
{7pt}% Space below
{}% Body font
{0pt}% Indent amount
{\bf}% Theorem head font
{:}% Punctuation after theorem head
{.5em}% Space after theorem head
{}% Theorem head spec (can be left empty, meaning `normal')

\newcommand\todo[1]{\textcolor{red}{($\bigstar$TODO$\bigstar$: #1)}}

\theoremstyle{amit}
\newtheorem{defn}{Definition}[section]
\newtheorem{prop}[defn]{Proposition}
\newtheorem{lem}[defn]{Lemma}
\newtheorem{thm}[defn]{Theorem}

\newtheorem{rmk}[defn]{Remark}
\newtheorem{ex}[defn]{Example}

\newcommand{\tatum}[1]           {{ \textcolor{olive}{[[ Tatum: {#1}]]}}}

\title{Poincar\'e Duality for Generalized Persistence Diagrams
of (co)Filtrations
\thanks {This work is partially funded by the 
Leverhulme Trust grant VP2-2021-008 awarded to the first author.}}
\author[1]{Amit Patel}
\author[1]{Tatum Rask}
\affil[1]{Department of Mathematics, Colorado State University}
\date{}

\begin{document}

\maketitle

\begin{abstract}
We dualize previous work on generalized persistence diagrams for filtrations to cofiltrations.
When the underlying space is a manifold, we express this duality as a Poincar\'e duality between their generalized persistence diagrams.
A heavy emphasis is placed on the recent discovery of functoriality of the generalized persistence diagram and its connection to Rota's Galois Connection Theorem.
\end{abstract}

%%%%%%%%%%%%%%%%%%%%%%%%%%%%%%%%%%%%%%%%%%%%%%
%%%%%%%%%%%%%%%%%%%%%%%%%%%%%%%%%%%%%%%%%%%%%%
\section{Introduction}
\label{sec:introduction}
%%%%%%%%%%%%%%%%%%%%%%%%%%%%%%%%%%%%%%%%%%%%%%
%%%%%%%%%%%%%%%%%%%%%%%%%%%%%%%%%%%%%%%%%%%%%%
Consider a filtration of a finite simplicial complex indexed by any finite poset.
For example, in classical persistence, the poset is totally ordered, and in multiparameter persistence, the poset is a product of totally ordered posets.
Its generalized persistent homology diagram, as defined in~\cite{McCPa20-edit-distance}, is an algebraic-combinatorial object describing the birth and death of cycles along the filtration.
In this paper, we begin a study of duality for generalized persistence diagrams by asking and answering the following question.
Is there a dual notion to a filtration that gives rise to a generalized persistent cohomology diagram in the same way?
The answer is yes, and we call it a \emph{cofiltration}.
When the underlying space is a manifold, we find a Poincar\'e duality between (co)filtrations and their generalized persistent (co)homology diagrams.
Similar questions have been asked of classical persistence diagrams~\cite{10.1145/2261250.2261287, deSilva_2011, extended_duality, sarah_duality}.

%%%%%%%%%%%%%%%%%%%%%%%%%%%%%%%%%%%%%%%%%%%%%%
\subsection{Prior Work}
%%%%%%%%%%%%%%%%%%%%%%%%%%%%%%%%%%%%%%%%%%%%%%
The inclusion-exclusion definition of the classical persistence diagram \cite{ELZ02, Robins99, 10.1117/12.279674} is a special case of the M\"obius inversion formula \cite{Patel2018}.
This discovery allows us to leverage the vast literature on the M\"obius inversion to generalize and uncover previously unknown structure~\cite{arxiv.2211.16642, GradedPersistenceDiagrams, McCPa20,McCPa20-edit-distance, GulenMcCleary, fasy_patel, kim2018generalized, Kim_Moore, virtual_dgm}.
We are no longer confined to totally ordered filtrations \cite{McCPa20-edit-distance, kim2018generalized, GulenMcCleary} nor to field coefficients~\cite{McCPa20}.
However, there is no consensus on how one should setup the M\"obius inversion when defining the (generalized) persistence diagram, e.g., compare \cite{GulenMcCleary} with~\cite{kim2018generalized}.
We believe questions of functoriality and stability will lead to a consensus. 

Kim and M\'emoli were the first to apply the M\"obius inversion to multiparameter persistence modules \cite{kim2018generalized}.
They define the persistence diagram as the M\"obius inversion of the generalized rank invariant.
Unlike classical persistence, these persistence diagrams have negative multiplicities.
%This makes the question of stability very difficult as we can not use existing methods for constructing a matching.

Functoriality of the persistence diagram is a long sought-after goal.
The first breakthrough came in \cite{McCPa20-edit-distance}.
Here, a persistence diagram is defined for any filtration of a finite simplicial complex indexed over any finite lattice.
This persistence diagram is defined as the M\"obius inversion of its birth-death function as opposed to its generalized rank invariant.
The idea of the birth-death function originates in the work of Henselman-Petrusek and Ghrist~\cite{GregH1, GregH2}.
The birth-death function behaves well with bounded-lattice functions leading to functoriality.
Again, this approach also produces negative multiplicities.
%Stability, in the form of an edit  distance, becomes an easy consequence of functoriality.

A further, very recent, breakthrough comes from G\"ulen and McCleary~\cite{GulenMcCleary}.
They discovered that functoriality, as achieved in \cite{McCPa20-edit-distance}, is equivalent to Rota's Galois Connection Theorem~\cite[Theorem 1]{Rot64}.
They go on to develop a functorial theory of persistence diagrams for persistence modules indexed by any finite poset again using the birth-death function.
They observe that the Galois connection offers a vast generalization of the interleaving distance, which works well with totally ordered sets, to arbitrary posets.

%%%%%%%%%%%%%%%%%%%%%%%%%%%%%%%%%%%%%%%%%%%%%%
\subsection{Contributions and Outline}
%%%%%%%%%%%%%%%%%%%%%%%%%%%%%%%%%%%%%%%%%%%%%%
The main contribution of this paper is an answer to the question of duality posed above.
In the process of answering this question, we provide a synthesis of ideas found in \cite{McCPa20-edit-distance} and \cite{GulenMcCleary}.
This synthesis emphasizes Rota's Galois Connection Theorem (Theorem~\ref{thm:rota}), and we hope this will be useful independent of duality.

We start (Section~\ref{sec:filtrations}) with an update to the persistence diagram for filtrations, as found in~\cite{McCPa20-edit-distance}, to the most general setting given what we know from \cite{GulenMcCleary}.
We then introduce and develop the dual theory for cofiltrations (Section~\ref{sec:cofiltrations}).
The following section (Section~\ref{sec:modules}) summaries the persistence diagram for persistence modules as found in~\cite{GulenMcCleary}.
By this point, there are two ways to construct a persistence diagram for a given (co)filtration. 
We show that both constructions give equivalent persistence diagrams (Section~\ref{sec:equivalence}).
We end with Poincar\'e duality (Section~\ref{sec:duality}).

%%%%%%%%%%%%%%%%%%%%%%%%%%%%%%%%%%%%%%%%%%%%%%
\subsection{Acknowledgements}
%%%%%%%%%%%%%%%%%%%%%%%%%%%%%%%%%%%%%%%%%%%%%%
We thank Alex McCleary for providing the proof for Proposition~\ref{prop:filtration_equivalence}.
The first author thanks Primoz Skraba for hosting him at Queen Mary University London (2022--2023) and for the many insightful discussions on the M\"obius inversion and persistent homology.

%%%%%%%%%%%%%%%%%%%%%%%%%%%%%%%%%%%%%%%%%%%%%%
%%%%%%%%%%%%%%%%%%%%%%%%%%%%%%%%%%%%%%%%%%%%%%
\section{Background}
\label{sec:background}
%%%%%%%%%%%%%%%%%%%%%%%%%%%%%%%%%%%%%%%%%%%%%%
%%%%%%%%%%%%%%%%%%%%%%%%%%%%%%%%%%%%%%%%%%%%%%
Before we start, we introduce the basic definitions and theorems at play in this paper. 

%%%%%%%%%%%%%%%%%%%%%%%%%%%%%%%%%%%%%%%%%%%%%%
\subsection{Posets}
%%%%%%%%%%%%%%%%%%%%%%%%%%%%%%%%%%%%%%%%%%%%%%
A \emph{poset} $P$ is a set with a reflexive, antisymmetric, and transitive
relation $\leq$.
A \emph{monotone function} from a poset $P$ to a poset $Q$ is a function $f : P \to Q$
such that if $a \leq b$ in $P$, then $f(a) \leq f(b)$ in $Q$.
A \emph{(monotone) Galois connection} from a poset $P$ to a poset $Q$, written
$f : P \leftrightarrows Q : g$, are monotone functions $f : P \to Q$ and $g : Q \to P$
such that $f(a) \leq x$ iff $a \leq g(x).$
Given two Galois connections $f : P \leftrightarrows Q : g$
and $h : Q \leftrightarrows R : i$, the two compositions form a Galois connection
$h \circ f : P \leftrightarrows R : g \circ i$.
Denote by $\Gal$ the category of posets and Galois connections.

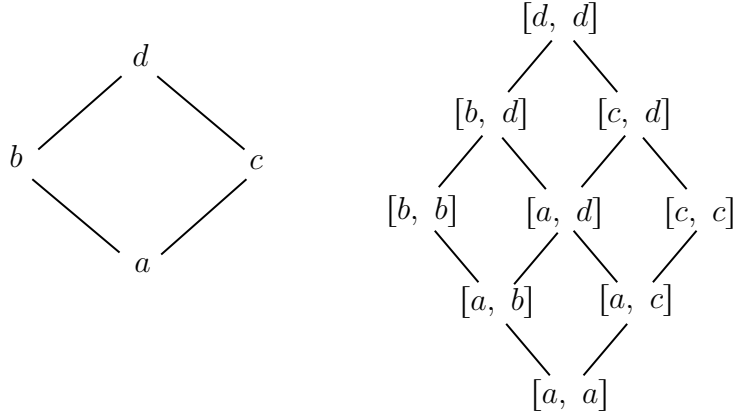
\begin{figure}
\begin{center}
% !TEX root = main.tex

\tikzset{every picture/.style={line width=0.75pt}} %set default line width to 0.75pt        

\begin{tikzpicture}[x=0.75pt,y=0.75pt,yscale=-1,xscale=1]

%Straight Lines [id:da5445702719803824] 
\draw    (223,96) -- (267,133) ;
%Straight Lines [id:da17860364428210018] 
\draw    (160,148) -- (204,185) ;
%Straight Lines [id:da2587685181684587] 
\draw    (268,148) -- (225,186) ;
%Straight Lines [id:da30470413900038285] 
\draw    (205,96) -- (163,133) ;
%Straight Lines [id:da7860604778984581] 
\draw    (460,221) -- (438,247) ;
%Straight Lines [id:da4774307279126826] 
\draw    (421,247) -- (399,221) ;
%Straight Lines [id:da5625936405574232] 
\draw    (385,200) -- (363,174) ;
%Straight Lines [id:da971211063528598] 
\draw    (455,200) -- (433,174) ;
%Straight Lines [id:da6074689907428141] 
\draw    (419,153) -- (397,127) ;
%Straight Lines [id:da22657469849682332] 
\draw    (490,153) -- (468,127) ;
%Straight Lines [id:da8745699683322208] 
\draw    (455,104) -- (433,78) ;
%Straight Lines [id:da3537447516348484] 
\draw    (495,174) -- (473,200) ;
%Straight Lines [id:da30311811175673076] 
\draw    (425,175) -- (403,201) ;
%Straight Lines [id:da8555164854208754] 
\draw    (387,126) -- (365,152) ;
%Straight Lines [id:da1276463573174087] 
\draw    (459,126) -- (437,152) ;
%Straight Lines [id:da552765290867244] 
\draw    (422,78) -- (400,104) ;

% Text Node
\draw (410,247.4) node [anchor=north west][inner sep=0.75pt]    {$[ a,\ a]$};
% Text Node
\draw (374,200.4) node [anchor=north west][inner sep=0.75pt]    {$[ a,\ b]$};
% Text Node
\draw (444,199.4) node [anchor=north west][inner sep=0.75pt]    {$[ a,\ c]$};
% Text Node
\draw (407,155.4) node [anchor=north west][inner sep=0.75pt]    {$[ a,\ d]$};
% Text Node
\draw (337,154.4) node [anchor=north west][inner sep=0.75pt]    {$[ b,\ b]$};
% Text Node
\draw (477,155.4) node [anchor=north west][inner sep=0.75pt]    {$[ c,\ c]$};
% Text Node
\draw (371,105.4) node [anchor=north west][inner sep=0.75pt]    {$[ b,\ d]$};
% Text Node
\draw (443,105.4) node [anchor=north west][inner sep=0.75pt]    {$[ c,\ d]$};
% Text Node
\draw (405,57.4) node [anchor=north west][inner sep=0.75pt]    {$[ d,\ d]$};
% Text Node
\draw (210,186.4) node [anchor=north west][inner sep=0.75pt]    {$a$};
% Text Node
\draw (147,129.4) node [anchor=north west][inner sep=0.75pt]    {$b$};
% Text Node
\draw (268,134.4) node [anchor=north west][inner sep=0.75pt]    {$c$};
% Text Node
\draw (209,78.4) node [anchor=north west][inner sep=0.75pt]    {$d$};

\end{tikzpicture}
\caption{Hasse diagrams of a lattice $P$ (left) and its corresponding interval 
lattice $\Int{P}$ (right)}
\label{fig:int-lattice}
\end{center}
\end{figure}

The \emph{interval} defined by a pair of elements $a \leq c$ in $P$ 
is the set $[a,c] := \{ b \in P : a \leq b \leq c  \}$. Denote by $\Int P$ the poset of all intervals in $P$. The product ordering is the natural ordering on $\Int P$.
That is, $[a,b] \leq [c,d]$ if $a \leq c$ and $b \leq d$.
The \emph{diagonal} of $\Int P$ is the subset
of all intervals of the form $[a,a]$.
See Figure~\ref{fig:int-lattice}.
A monotone function $f : P \to Q$ induces a monotone function $\Int f : \Int P \to \Int Q$ where $\Int f [a,b] := \big[ f(a), f(b) \big]$.

The following lemma implies that the operation of taking intervals
is an endofunctor $\Int : \Gal \to \Gal$.

\begin{lem}
    If $f : P \leftrightarrows Q : g$ is a Galois connection,
    then $\Int f : \Int P \leftrightarrows \Int Q : \Int g$  is a Galois
    connection.
\end{lem}
\begin{proof}
The statement follows easily from the product ordering on intervals:
    \begin{align*}
        \Int f [a,b] \leq [x,y] & \Leftrightarrow  [a,b] \leq \Int g [x,y] \\
        f(a) \leq x \text{ and } f(b) \leq y & \Leftrightarrow  a \leq g(x) \text{ and }
        b \leq g(y).
    \end{align*}
\end{proof}

\begin{rmk}
For historical reasons \cite{stability-persistence}, it is common to impose
the containment ordering on $\Int P$ \cite{kim2018generalized}.
That is, $[a,b] \leq [c,d]$ if $c \leq a$ and $b \leq d$.
The containment ordering is the right ordering when
taking the M\"obius inversion of the generalized rank invariant.
Unfortunately, this ordering 
fails to produce an endofunctor on $\Gal$ making unavailable Theorem~\ref{thm:rota}.
That is, the M\"obius inversion of the
generalized rank invariant is not functorial
with respect to Galois connections between (co)filtrations; see Propositions~\ref{prop:DGMH}, \ref{prop:DGMC}, and~\ref{prop:dgm_galois}.
\end{rmk}

%%%%%%%%%%%%%%%%%%%%%%%%%%%%%%%%%%%%%%%%%%%%%%
\subsection{M\"obius Inversion}
%%%%%%%%%%%%%%%%%%%%%%%%%%%%%%%%%%%%%%%%%%%%%%
The M\"obius inversion has its roots in number theory, but the modern theory was introduced by Rota \cite{Rot64}.
The literature is rich with deep ideas waiting to find their way into persistent homology.
See, for example, \cite{CRAPO1966126, BACLAWSKI1975125, rota_euler}.
This paper focuses on one such idea: Rota's Galois Connection Theorem.
As we shall see, it translates to functoriality of the persistence diagram.

Let $P$ be a finite poset and $m : P \to \Z$ an integral function.
Then, there is a unique function $\partial m : P \to \Z$, called the
\emph{M\"obius inversion} of $m$, such that for all $b \in P$,
    $$m(b) = \sum_{a : a \leq b} \partial m(a).$$
The M\"obius inversion $\partial m$ can be interpreted as the combinatorial
derivative of $m$, hence the notation.

Fix a Galois connection $f : P \leftrightarrows Q : g$.
The \emph{pullback} of a function $m : P \to \Z$ along the Galois connection 
is the function $g^\sharp : Q \to \Z$ defined
as $g^\sharp = m \circ g$.
The \emph{pushforward} of a function $m : P \to \Z$ along the Galois connection 
is the function $f_\sharp : Q \to \Z$ defined as
$$f_\sharp (x) = \sum_{a \in f^{-1}(x)} m(a).$$
If $f^{-1}(x)$ is empty, then we interpret the sum as~$0$.

\begin{thm}[Rota's Galois Connection Theorem~\cite{GulenMcCleary}]
\label{thm:rota}
Let $m : P \to \Z$ and $n : Q \to \Z$ be two functions and $f : P \leftrightarrows Q :  g$
a Galois connection. If $n = g^\sharp m$, then $\partial n = f_\sharp (\partial m)$.
\end{thm}

The theorem implies a functor between the following two categories:
    \begin{itemize}
    \item Denote by $\mathbb{A}$ the category of functions $m : P \to \Z$, over all finite
    posets.
    A morphism from $m  : P \to \Z$ to $n : Q \to \Z$ is a Galois 
    connection $f  : P \leftrightarrows Q :  g$ such that $n = g^\sharp m$.
    We think of $\mathbb{A}$ as the category of birth-death functions described 
    in three different ways later in this paper.

    \item Denote by $\mathbb{B}$ the category of functions $m : P \to \Z$, over all finite
    posets.
    A morphism from $m  : P \to \Z$ to $n : Q \to \Z$ is a Galois 
    connection $f  : P \leftrightarrows Q :  g$ such that $n = f_\sharp m$.
    We think of $\mathbb{B}$ as the category of persistence diagrams described 
    in three different ways later in this paper.
    \end{itemize}
The theorem says that the M\"obius inversion maps a morphism in $\mathbb{A}$ to
a morphism in $\mathbb{B}$.

%%%%%%%%%%%%%%%%%%%%%%%%%%%%%%%%%%%%%%%%%%%%%%
\subsection{Cohomology with Compact Support}
%%%%%%%%%%%%%%%%%%%%%%%%%%%%%%%%%%%%%%%%%%%%%%

We now review cohomology with compact support and state a version of Poincar\'e duality~\cite[pg.\ 242]{Hatcher}.
Below, $\Vect$ is the category of $\field$-vector spaces, for some fixed field  $\field$.

Fix a topological space $X$ and consider its
poset $\mathcal{K}(X)$ of all compact
subsets ordered by inclusion.
We now define a functor $\mathcal{H}^\ast : \mathcal{K}(X) \to \Vect$.
For a compact set $A$, let $\mathcal{H}^\ast(A)$ be the relative singular cohomology group $H^\ast(X, X \setminus A)$.
For $A \subseteq B$, the inclusion of pairs $(X, X \setminus B) \hookrightarrow (X, X \setminus A)$ induces a map $H^\ast(X, X \setminus A) \to H^\ast(X, X \setminus B)$.
\emph{Cohomology with compact support} is the colimit of this functor.
That is, $H^\ast_\com(X) := \mathsf{colim\,} \mathcal{H}^\ast$.

Now suppose $X$ is an open subset of a larger space $Y$.
Every compact subset of~$X$ is also a compact subset of $Y$.
That is, $\mathcal{K}(X) \hookrightarrow \mathcal{K}(Y)$.
Every compact set $A \subseteq X$ is contained in an open set $U \subseteq X$.
The following maps induced by inclusion are, by excision, isomorphisms:
    \begin{equation*}
    \begin{tikzcd}
    H^\ast(X, X \setminus A) \ar[rr, "\cong"] && H^\ast(U, U \setminus A) 
    &&  H^\ast(Y, Y \setminus A). \ar[ll, "\cong"']
    \end{tikzcd}
    \end{equation*}
Thus, the universal property of the colimit $H^\ast_\com(X)$ induces a canonical map
$H^\ast_\com(X) \to H^\ast_\com(Y)$.

\begin{thm}[Poincar\'e  Duality]
\label{thm:poincare}
Let $Y$  be an $\field$-orientable, not necessarily compact, $m$-manifold.
Then, there is a canonical isomorphism $H^i_\com(Y) \cong H_{m-i}(Y)$.
Further, for any open subset $X \subseteq Y$, the following diagram, where the horizontal maps are induced by inclusion, commutes:
    \begin{equation*}
    \begin{tikzcd}
	H^i_\com(X) \ar[d, "\cong"]  \ar[rr] && H^i_\com (Y) \ar[d,  "\cong"] \\
	H_{m-i}(X) \ar[rr] && H_{m-i}(Y).
    \end{tikzcd}
    \end{equation*}
    \end{thm}

%%%%%%%%%%%%%%%%%%%%%%%%%%%%%%%%%%%%%%%%%%%%%%
%%%%%%%%%%%%%%%%%%%%%%%%%%%%%%%%%%%%%%%%%%%%%%
\section{Filtrations}
\label{sec:filtrations}
%%%%%%%%%%%%%%%%%%%%%%%%%%%%%%%%%%%%%%%%%%%%%%
%%%%%%%%%%%%%%%%%%%%%%%%%%%%%%%%%%%%%%%%%%%%%%
Fix a finite simplicial complex $K$ and denote by $\Delta K$ its poset
of subcomplexes ordered by inclusion.
Recall a \emph{subcomplex} is a subset of simplices $A \subseteq K$ such that if 
$\tau \in A$ and $\sigma \leq \tau$, then $\sigma \in A$.

A \emph{filtration} of $K$ indexed by a poset $P$ is a monotone function $F : P \to \Delta K$.
See Figure~\ref{fig:filtration}.
A \emph{filtration-preserving} morphism from a filtration 
$F : P \to \Delta K$ to a filtration $G : Q \to \Delta K$
is a Galois connection $f : P \leftrightarrows Q : g$ such that
$G = F \circ g$.
That is, $G$ is the pullback of $F$ along the Galois connection.
We write such a morphism as $f : F \leftrightarrows G : g$.

%%%%%%%%%%%%%%%%%%%%%%%%%%%%%%%%%%%%%%%%%%%%%%
\subsection{Birth-Death Function}
%%%%%%%%%%%%%%%%%%%%%%%%%%%%%%%%%%%%%%%%%%%%%%
Fix a filtration $F : P \to \Delta K$.
For every $a \in P$, let $C_\bullet F(a)$ be the usual simplicial chain complex
over coefficients in a fixed field $\field$.
For every $a \leq b$, we have an inclusion of chain complexes:
    \begin{equation*}
    \begin{tikzcd}
    \cdots \ar[r] & C_{d+1} F(a) \ar[r, "\partial_{d+1}"] 
    \ar[d, hookrightarrow ] &  C_{d} F(a) \ar[r, "\partial_{d}"] 
    \ar[d, hookrightarrow ] & 
    \ar[r] 
    \ar[d, hookrightarrow ] C_{d-1} F(a) & \cdots \\
    \cdots \ar[r] & C_{d+1} F(b) \ar[r, "\partial_{d+1}"] &  
    C_{d} F(b) \ar[r, "\partial_{d}"] & C_{d-1} F(b) \ar[r] & \cdots .
    \end{tikzcd}
    \end{equation*}
Denote by $Z_d F(a)$ and $B_d F(a)$
the $d$-th cycle space and $d$-th boundary space, respectively, at~$a$.
Cycles and boundaries include along the filtration.
That is, for $a \leq b$, $Z_d F(a) \subseteq Z_d F(b)$ and 
$B_d F(a) \subseteq B_d F(b)$.

\begin{defn}
For $a \leq b$, let $Z_d F(a) \cap B_d F(b)$ be the pullback of the following
diagram of $\field$-vector spaces:
    \begin{equation*}
    \begin{tikzcd}
    Z_d F(a) \cap B_d F(b) \ar[rr, dashrightarrow] \ar[d, dashrightarrow]
    && Z_d F(a) \ar[d, hookrightarrow] \\
    B_d F(b) \ar[rr, hookrightarrow] 
    && Z_d F(b).
    \end{tikzcd}
    \end{equation*}
The \define{$d$-th birth-death function} $\BD_d F :  \Int P \to \Z$ is defined as  
    $$\BD_d F [a,b] = \dim \big( Z_d F(a) \cap B_d F(b) \big).$$
\end{defn}

\begin{figure}
\begin{center}
\include{figures-tex/filtration.tex}
\caption{A filtration $F: P \to \Delta K$ (left), its corresponding birth-death function $\BD_1  F: \Int{P} \to \Z$ (middle), and the persistent homology diagram $\partial \BD_1 F : \Int{P} \to \Z$ (right)}
\label{fig:filtration}
\end{center}
\end{figure}
See Figure~\ref{fig:filtration} for an example of a birth-death function. 

    \begin{prop}\label{prop:BDH}
    If $f : F \leftrightarrows G: g$ is a filtration-preserving morphism,
    then $\BD_\ast G = (\Int g)^\sharp \big( \BD_\ast F \big)$.
    \end{prop}
    \begin{proof}
    We want to show $\BD_\ast G [x,y] = 
    \BD_\ast F \big[ g(x), g(y) \big]$.
    This follows immediately from the fact that
    $G(x) = F\big( g(x) \big)$ and
    $G(y) = F \big ( g(y) \big)$.
    In other words, the two pullback diagrams defining $\BD_\ast G [x,y]$ and $\BD_\ast F \big[ g(x), g(y) \big]$ are the same.
    \end{proof}

%%%%%%%%%%%%%%%%%%%%%%%%%%%%%%%%%%%%%%%%%%%%%%
\subsection{Persistent Homology Diagram}
%%%%%%%%%%%%%%%%%%%%%%%%%%%%%%%%%%%%%%%%%%%%%%

We interpret $\BD_d F[a,b]$ as cycles born by $a$ and dead by $b$.
However, the persistence diagram of $F$ should capture cycles born at $a$
and dead at $b$.
We achieve this by taking the combinatorial derivative.

\begin{defn}
    \label{defn:cycles_boundaries}
    The \define{$d$-th persistent homology diagram} of a filtration $F : P \to \Delta K$
    is the M\"obius inversion $\partial \BD_d F : \Int P \to \Z$ of its $d$-th
    birth-death function.
\end{defn}

When $P$ is totally ordered, $\partial \BD_d F$ is the classical persistent homology
diagram \cite[\S 9]{McCPa20-edit-distance}.
See Figure~\ref{fig:filtration} for an example of a persistent homology diagram.
Notice the negative multiplicities.

    \begin{prop}\label{prop:DGMH}
    If $f : F \leftrightarrows G: g$ is a filtration-preserving morphism,
    then $\partial \BD_d G = (\Int f)_\sharp \big( \partial \BD_d F \big)$.
    \end{prop}
    \begin{proof}
    The statement follows from Proposition~\ref{prop:BDH} and Theorem~\ref{thm:rota}.
    \end{proof}

%%%%%%%%%%%%%%%%%%%%%%%%%%%%%%%%%%%%%%%%%%%%%%
%%%%%%%%%%%%%%%%%%%%%%%%%%%%%%%%%%%%%%%%%%%%%%
\section{Cofiltrations}
\label{sec:cofiltrations}
%%%%%%%%%%%%%%%%%%%%%%%%%%%%%%%%%%%%%%%%%%%%%%
%%%%%%%%%%%%%%%%%%%%%%%%%%%%%%%%%%%%%%%%%%%%%%
Fix a finite simplicial complex $K$.
A \emph{supcomplex} of $K$ is a subset of simplices $A \subseteq K$ such that if $\sigma \in A$ and $\sigma \leq \tau$, then $\tau \in A$.
Denote by $\nabla K$ the poset
of supcomplexes ordered by inclusion.

A \emph{cofiltration} of $K$ indexed by a poset $P$ is a monotone function
$F : P \to \nabla K$.
See Figure~\ref{fig:mob-ph-ex}.
A \emph{cofiltration-preserving} morphism from a cofiltration 
$F : P \to \nabla K$ to a cofiltration $G : Q \to \nabla K$
is a Galois connection $f : P \leftrightarrows Q : g$ such that
$G = F \circ g$.
That is, $G$ is the pullback of $F$ along the Galois connection.
We write such a morphism as $f : F \leftrightarrows G : g$.

%%%%%%%%%%%%%%%%%%%%%%%%%%%%%%%%%%%%%%%%%%%%%%
\subsection{Birth-Death Function}
%%%%%%%%%%%%%%%%%%%%%%%%%%%%%%%%%%%%%%%%%%%%%%
Fix a cofiltration $F : P \to \nabla K$.
For every $a \in P$, consider the chain complex $C_\bullet F(a)$ as defined in the previous section but with one caveat.
A simplex $\tau \in F(a)$ may not have all its faces in $F(a)$.
This means that if $\sigma$ is on the boundary of $\tau$ and $\sigma$ is not in $F(a)$, then the boundary operator applied to $\tau$ must send $\sigma$ to $0$.

For every $a \leq b$, we have a surjection of chain complexes as follows:
    \begin{equation*}
    \begin{tikzcd}
    \cdots \ar[r] & C_{d+1} F(b) \ar[r, "\partial_{d+1}"] 
    \ar[d, twoheadrightarrow, "i_{d+1}"] &  C_{d} F(b) \ar[r, "\partial_{d}"] 
    \ar[d, twoheadrightarrow, "i_d"] & 
    \ar[r] 
    \ar[d, twoheadrightarrow, "i_{d-1}"] C_{d-1} F(b) & \cdots \\
    \cdots \ar[r] & C_{d+1} F(a) \ar[r, "\partial_{d+1}"] &  
    C_{d} F(a) \ar[r, "\partial_{d}"] & C_{d-1} F(a) \ar[r] & \cdots .
    \end{tikzcd}
    \end{equation*}
The map $i_d$ is generated by sending every
$d$-simplex $\tau \in F(b)$ to itself, if it exists in $F(a)$, or to~$0$, otherwise.
This makes every $i_d$ surjective.
Next, we want to show $i_{d-1} \circ \partial_d = \partial_d \circ i_d$.
Suppose $\sigma$ is on the boundary of $\tau \in F(b)$ and $\sigma \in F(b)$.
If $\sigma \in F(a)$, then $\tau \in F(a)$ because~$F(a)$ is a supcomplex, and so the square commutes.
If $\sigma \notin F(a)$, then $i_{d-1} \circ \partial_d$ and $\partial_d \circ i_d$ applied to $\tau$ sends $\sigma$ to~$0$, and so the square commutes.

Now dualize everything above.
That is, let $C^\bullet F(a)$ be the cochain complex
$\Hom \big( C_\bullet, \field \big)$.
For $a \leq b$, the surjective chain map $i_\bullet$ above
dualizes to an injective cochain map $j^\bullet$:
    \begin{equation*}
    \begin{tikzcd}
    \cdots \ar[r] & C^{d-1} F(a) \ar[r, "\delta^{d-1}"] 
    \ar[d, hookrightarrow, "j^{d-1}"] &  C^{d} F(a) \ar[r, "\delta^{d}"] 
    \ar[d, hookrightarrow, "j^d"] & \ar[r] 
    \ar[d, hookrightarrow, "j_{d+1}"] C^{d+1} F(a) & \cdots \\
    \cdots \ar[r] & C^{d-1} F(b) \ar[r, "\delta^{d-1}"] &  
    C^{d} F(b) \ar[r, "\delta^{d}"] & C^{d+1} F(b) \ar[r] & \cdots.
    \end{tikzcd}
    \end{equation*}
Denote by $Z^d F(a)$ and $B^d F(a)$
the $d$-th cocycle space and $d$-th coboundary space, respectively, at $a$.
The cohomology theory we have just described is \emph{cellular cohomology with compact support} \cite{shepard}.

\begin{ex}
Suppose $K$ is the $1$-simplex and $F : P \to \nabla K$ a two-step cofiltration where
$F(a)$ is the $1$-simplex without its two vertices and $F(b)$ is the $1$-simplex with both its vertices.
Then, $H^1_\com F(a) := \sfrac{Z^1 F(a)}{B^1 F(a)}  \cong \field$ wheres
the singular cohomology $H^1 \big( | F(a) | \big)$ of the underlying space,  which is homeomorphic to the real line, is $0$.
\end{ex}

\begin{figure}
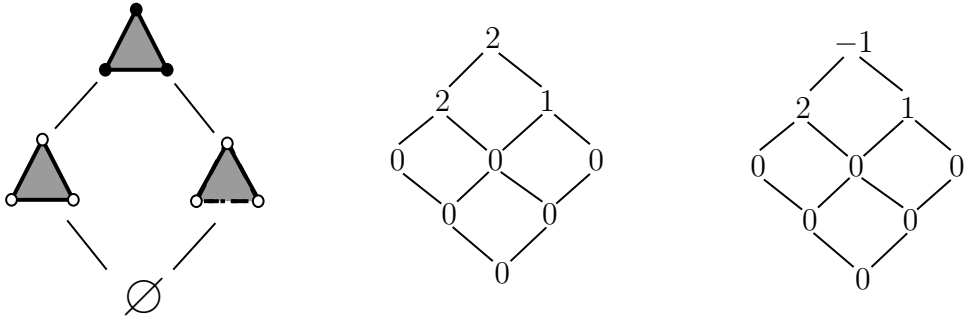

\begin{center}
\include{figures-tex/mob-ph-ex.tex}
\caption{A cofiltration $F: P \to \nabla K$ (left), its corresponding birth-death function $\BD^1 F: \Int P \to \mathbb{Z}$ (middle), and the  persistent cohomology diagram $\partial \BD^1 F: \Int P \to \mathbb{Z}$ (right)}
\label{fig:mob-ph-ex}
\end{center}
\end{figure}

\begin{defn}
For $a \leq b$, let $Z^d F(a) \cap B^d F(b)$ be the pullback of the following
diagram of $\field$-vector spaces:
    \begin{equation*}
    \begin{tikzcd}
    Z^d F(a) \cap B^d F(b) \ar[rr, dashrightarrow] \ar[d, dashrightarrow]
    && Z^d F(a) \ar[d, hookrightarrow] \\
    B^d F(b) \ar[rr, hookrightarrow] 
    && Z^d F(b).
    \end{tikzcd}
    \end{equation*}
The \define{$d$-th birth-death function} $\BD^d F :  \Int P \to \Z$ is defined as  
    $$\BD^d F [a,b] = \dim \big( Z^d F(a) \cap B^d F(b) \big).$$
\end{defn}

See Figure~\ref{fig:mob-ph-ex} for an example of a birth-death function.

    \begin{prop}\label{prop:BDC}
    If $f : F \leftrightarrows G: g$ is a cofiltration-preserving morphism,
    then $\BD^\ast G = (\Int g)^\sharp \big( \BD^\ast F \big)$.
    \end{prop}
    \begin{proof}
    We want to show $\BD^\ast G [x,y] = 
    \BD^\ast F \big[ g(x), g(y) \big]$.
    This follows immediately from the fact that
    $G(x) = F\big( g(x) \big)$ and
    $G(y) = F \big ( g(y) \big)$.
    In other words, the pullback diagrams defining $\BD^\ast G [x,y]$ and $\BD^\ast F \big[ g(x), g(y) \big]$ are the same.
    \end{proof}

%%%%%%%%%%%%%%%%%%%%%%%%%%%%%%%%%%%%%%%%%%%%%%
\subsection{Persistent Cohomology Diagram}
%%%%%%%%%%%%%%%%%%%%%%%%%%%%%%%%%%%%%%%%%%%%%%
We interpret $\BD^d F[a,b]$ as cocycles born by $a$ and dead by $b$.
However, the persistence diagram of $F$ should capture cocycles born at $a$
and dead at $b$.
We achieve this by taking the combinatorial derivative.

\begin{defn}
    \label{defn:cocyles_coboundaries}
    The \define{$d$-th persistent cohomology diagram} of a cofiltration $F : P \to \nabla K$
    is the M\"obius inversion $\partial \BD^d F : \Int P \to \Z$ of its $d$-th
    birth-death function.
\end{defn}

See Figure~\ref{fig:mob-ph-ex} for an example of a persistent cohomology diagram.

    \begin{prop} \label{prop:DGMC}
    If $f : F \leftrightarrows G: g$ is a cofiltration preserving morphism,
    then $\partial \BD^d G = (\Int f)_\sharp \big( \partial \BD^d F \big)$.
    \end{prop}
    \begin{proof}
    The statement follows from Proposition~\ref{prop:BDC} and Theorem~\ref{thm:rota}.
    \end{proof}

%%%%%%%%%%%%%%%%%%%%%%%%%%%%%%%%%%%%%%%%%%%%%%
%%%%%%%%%%%%%%%%%%%%%%%%%%%%%%%%%%%%%%%%%%%%%%
\section{Persistence Modules}
\label{sec:modules}
%%%%%%%%%%%%%%%%%%%%%%%%%%%%%%%%%%%%%%%%%%%%%%
%%%%%%%%%%%%%%%%%%%%%%%%%%%%%%%%%%%%%%%%%%%%%%
Here we provide a summary of ideas detailed in~\cite{GulenMcCleary} that will be useful later.
Below $\Vect$ is the category of finite dimensional $\field$-vector spaces and $\field$-linear maps.

A \emph{persistence module} indexed by a poset $P$, or a \emph{$P$-module}, is a functor $M : P \to \Vect$.
A \emph{module-preserving morphism} from a persistence module~$M : P \to \Vect$ to a persistence module~$N : Q \to \Vect$ is a Galois connection~$f : P \leftrightarrows Q : g$  such that $N = M  \circ f$.
That is,~$N$ is the pullback of $M$ along the Galois connection.
We write such a morphism as $f : M \leftrightarrows N : g$.

For $a \in P$, let $\field^{\uparrow a}$ be the persistence module
    \begin{equation*}
    \field^{\uparrow a}(b) =
    \begin{cases}
    \field \text{ if $a \leq b$} \\
    0 \text{ otherwise}
    \end{cases}
    \end{equation*}
and $\field^{\uparrow a}(b \leq c) = 1_{\field}$
for all $a \leq b \leq c$.
A $P$-module is \emph{free} if it is the direct sum of modules of the form $\field^{\uparrow a}$~\cite{miller_hom_algebra}.

%%%%%%%%%%%%%%%%%%%%%%%%%%%%%%%%%%%%%%%%%%%%%%
\subsection{Birth-Death Function}
%%%%%%%%%%%%%%%%%%%%%%%%%%%%%%%%%%%%%%%%%%%%%%

\begin{figure}
\begin{center}
% !TEX root = main.tex

\tikzset{every picture/.style={line width=0.75pt}} %set default line width to 0.75pt        

\begin{tikzpicture}[x=0.75pt,y=0.75pt,yscale=-1,xscale=1]
%uncomment if require: \path (0,243); %set diagram left start at 0, and has height of 243

%Straight Lines [id:da019345910176327807] 
\draw    (168,129) -- (200,94) ;
%Straight Lines [id:da373333037599632] 
\draw    (226,194) -- (259,156) ;
%Straight Lines [id:da027948797834237205] 
\draw    (257,129) -- (227,94) ;
%Straight Lines [id:da43313303876452736] 
\draw    (205,193) -- (170,157) ;
%Straight Lines [id:da8160900265582636] 
\draw [color={rgb, 255:red, 128; green, 128; blue, 128 }  ,draw opacity=1 ]   (269.02,79.64) -- (415.1,79.64) ;
\draw [shift={(417.1,79.64)}, rotate = 180] [color={rgb, 255:red, 128; green, 128; blue, 128 }  ,draw opacity=1 ][line width=0.75]    (10.93,-3.29) .. controls (6.95,-1.4) and (3.31,-0.3) .. (0,0) .. controls (3.31,0.3) and (6.95,1.4) .. (10.93,3.29)   ;
%Straight Lines [id:da012599674223931157] 
\draw [color={rgb, 255:red, 128; green, 128; blue, 128 }  ,draw opacity=1 ]   (269.02,202.23) -- (415.1,202.23) ;
\draw [shift={(417.1,202.23)}, rotate = 180] [color={rgb, 255:red, 128; green, 128; blue, 128 }  ,draw opacity=1 ][line width=0.75]    (10.93,-3.29) .. controls (6.95,-1.4) and (3.31,-0.3) .. (0,0) .. controls (3.31,0.3) and (6.95,1.4) .. (10.93,3.29)   ;
%Curve Lines [id:da6656176749738139] 
\draw [color={rgb, 255:red, 128; green, 128; blue, 128 }  ,draw opacity=1 ]   (267,155) .. controls (295.52,170.12) and (448.8,170.67) .. (490.19,153.78) ;
\draw [shift={(492,153)}, rotate = 155.33] [color={rgb, 255:red, 128; green, 128; blue, 128 }  ,draw opacity=1 ][line width=0.75]    (10.93,-3.29) .. controls (6.95,-1.4) and (3.31,-0.3) .. (0,0) .. controls (3.31,0.3) and (6.95,1.4) .. (10.93,3.29)   ;
%Curve Lines [id:da4498387400320052] 
\draw [color={rgb, 255:red, 128; green, 128; blue, 128 }  ,draw opacity=1 ]   (178,126) .. controls (212,116.19) and (363.32,115.59) .. (402.51,126.03) ;
\draw [shift={(404.22,126.51)}, rotate = 196.78] [color={rgb, 255:red, 128; green, 128; blue, 128 }  ,draw opacity=1 ][line width=0.75]    (10.93,-3.29) .. controls (6.95,-1.4) and (3.31,-0.3) .. (0,0) .. controls (3.31,0.3) and (6.95,1.4) .. (10.93,3.29)   ;
%Straight Lines [id:da10566270897346541] 
\draw    (471,193) -- (504,155) ;
%Straight Lines [id:da8163258477656907] 
\draw    (450,192) -- (415,156) ;
%Straight Lines [id:da38185801882012216] 
\draw    (414,129) -- (446,94) ;
%Straight Lines [id:da9147247666662717] 
\draw    (503,129) -- (473,94) ;

% Text Node
\draw (207.09,71.86) node [anchor=north west][inner sep=0.75pt]  [font=\normalsize]  {$\mathbb{F}^{3}$};
% Text Node
\draw (158.27,132.36) node [anchor=north west][inner sep=0.75pt]  [font=\normalsize]  {$\mathbb{F}^{2}$};
% Text Node
\draw (253.31,133.22) node [anchor=north west][inner sep=0.75pt]  [font=\normalsize]  {$\mathbb{F}$};
% Text Node
\draw (209.67,197.4) node [anchor=north west][inner sep=0.75pt]  [font=\normalsize]  {$0$};
% Text Node
\draw (454.09,74.86) node [anchor=north west][inner sep=0.75pt]  [font=\normalsize]  {$0$};
% Text Node
\draw (403.27,132.36) node [anchor=north west][inner sep=0.75pt]  [font=\normalsize]  {$\mathbb{F}^{2}$};
% Text Node
\draw (498.31,133.22) node [anchor=north west][inner sep=0.75pt]  [font=\normalsize]  {$\mathbb{F}$};
% Text Node
\draw (454.67,197.4) node [anchor=north west][inner sep=0.75pt]  [font=\normalsize]  {$0$};
% Text Node
\draw (330.16,62.82) node [anchor=north west][inner sep=0.75pt]  [font=\normalsize,color={rgb, 255:red, 128; green, 128; blue, 128 }  ,opacity=1 ]  {$0$};
% Text Node
\draw (330.16,184.87) node [anchor=north west][inner sep=0.75pt]  [font=\normalsize,color={rgb, 255:red, 128; green, 128; blue, 128 }  ,opacity=1 ]  {$0$};
% Text Node
\draw (370.44,149.45) node [anchor=north west][inner sep=0.75pt]  [font=\normalsize,color={rgb, 255:red, 128; green, 128; blue, 128 }  ,opacity=1 ]  {$\textrm{id}$};
% Text Node
\draw (274.28,97.34) node [anchor=north west][inner sep=0.75pt]  [font=\normalsize,color={rgb, 255:red, 128; green, 128; blue, 128 }  ,opacity=1 ]  {$\textrm{id} \ $};
% Text Node
\draw (150,77.4) node [anchor=north west][inner sep=0.75pt]  [font=\tiny]  {$\begin{pmatrix}
1 & 0\\
0 & 1\\
0 & 0
\end{pmatrix}$};
% Text Node
\draw (245,76.4) node [anchor=north west][inner sep=0.75pt]  [font=\tiny]  {$\begin{pmatrix}
0\\
0\\
1
\end{pmatrix}$};

\end{tikzpicture}
\caption{A persistence module $M$ (right) %corresponding to the cofiltration $F$ in Figure \ref{fig:mob-ph-ex} 
and a free module $F$ (left). A free presentation of $M$ is given by the natural transformation $\phi : F \Rightarrow M$ shown in grey.}
\label{fig:presentation}
\end{center}
\end{figure}
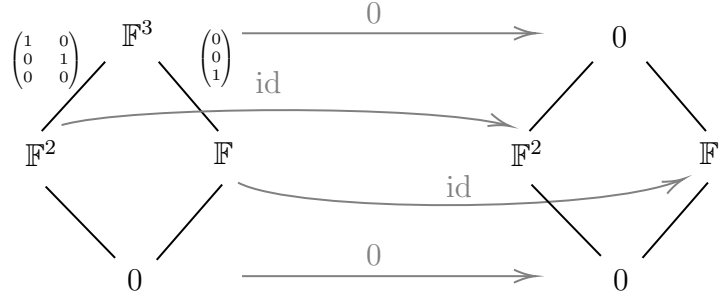

Fix a module $M : P \to \Vect$.
A \emph{presentation} of $M$ is a surjective natural transformation $\phi : F \Rightarrow M$ from a free $P$-module~$F$.
A presentation always exists, but it is not unique.
For example, one can introduce an arbitrary number of extra generators to $F$ as long as those generators are killed immediately by relations.
See Figure~\ref{fig:presentation} for an example of a presentation of a persistence module.

\begin{defn}
Fix a presentation $\phi : F \Rightarrow M$.
For $a \leq b$, let $F (a) \cap \ker \phi_b$ be the pullback of the following diagram of $\field$-vector spaces:
    \begin{equation*}
    \begin{tikzcd}
    F (a) \cap \ker \phi_b \ar[rr, dashrightarrow] \ar[d, dashrightarrow]
    && \ker \phi_b \ar[d, hookrightarrow] \\
    F(a)  \ar[rr, hookrightarrow] 
    && F (b).
    \end{tikzcd}
    \end{equation*}
The \define{birth-death} function $\BD \phi: \Int P \to \Z$ of a presentation $\phi : F \Rightarrow M$ is defined as   
$$\BD \phi [a,b]  = \dim \big( F (a) \cap \ker \phi_b \big).$$
\end{defn}

\begin{figure}
\begin{center}
% !TEX root = main.tex

\tikzset{every picture/.style={line width=0.75pt}} %set default line width to 0.75pt        

\begin{tikzpicture}[x=0.75pt,y=0.75pt,yscale=-1,xscale=1]
%uncomment if require: \path (0,222); %set diagram left start at 0, and has height of 222

%Straight Lines [id:da13342641045286885] 
\draw    (216,72) -- (198,88) ;
%Straight Lines [id:da9070629253433733] 
\draw    (190,103) -- (172.66,117.04) ;
%Straight Lines [id:da8819023574929068] 
\draw    (245,102) -- (225,120) ;
%Straight Lines [id:da9955419260855782] 
\draw    (218,131) -- (202,146) ;
%Straight Lines [id:da6277422605421026] 
\draw    (247,159) -- (229,176) ;
%Straight Lines [id:da030973595344556903] 
\draw    (269,132) -- (252,147) ;
%Straight Lines [id:da7977101110895766] 
\draw    (244.64,88.02) -- (224,72) ;
%Straight Lines [id:da08989475520577628] 
\draw    (220,176) -- (199.77,158.8) ;
%Straight Lines [id:da5624962247861756] 
\draw    (194,147) -- (174.77,131.8) ;
%Straight Lines [id:da6873909491187158] 
\draw    (217,120) -- (196.77,102.8) ;
%Straight Lines [id:da6905810315199814] 
\draw    (269,117) -- (251,102) ;
%Straight Lines [id:da755102022782717] 
\draw    (245,146) -- (224,131) ;
%Straight Lines [id:da28214667736453536] 
\draw    (428,72) -- (410,88) ;
%Straight Lines [id:da36021598964829127] 
\draw    (402,103) -- (384.66,117.04) ;
%Straight Lines [id:da7749751593966143] 
\draw    (457,102) -- (437,120) ;
%Straight Lines [id:da2696800079905475] 
\draw    (430,131) -- (414,146) ;
%Straight Lines [id:da5823989019065974] 
\draw    (459,159) -- (441,176) ;
%Straight Lines [id:da5036379476071395] 
\draw    (481,132) -- (464,147) ;
%Straight Lines [id:da7094343026977843] 
\draw    (456.64,88.02) -- (436,72) ;
%Straight Lines [id:da510483449777392] 
\draw    (432,176) -- (411.77,158.8) ;
%Straight Lines [id:da9281970088709215] 
\draw    (406,147) -- (386.77,131.8) ;
%Straight Lines [id:da38671684908067805] 
\draw    (429,120) -- (408.77,102.8) ;
%Straight Lines [id:da0980984673316041] 
\draw    (481,117) -- (463,102) ;
%Straight Lines [id:da6262441704506005] 
\draw    (457,146) -- (436,131) ;

% Text Node
\draw (431.33,175.08) node [anchor=north west][inner sep=0.75pt]    {$0$};
% Text Node
\draw (404.6,145.37) node [anchor=north west][inner sep=0.75pt]    {$0$};
% Text Node
\draw (455.33,145.37) node [anchor=north west][inner sep=0.75pt]    {$0$};
% Text Node
\draw (427.9,118.43) node [anchor=north west][inner sep=0.75pt]    {$0$};
% Text Node
\draw (378.55,118.43) node [anchor=north west][inner sep=0.75pt]    {$0$};
% Text Node
\draw (478.63,118.43) node [anchor=north west][inner sep=0.75pt]    {$0$};
% Text Node
\draw (401.17,88.72) node [anchor=north west][inner sep=0.75pt]    {$2$};
% Text Node
\draw (453.95,88.72) node [anchor=north west][inner sep=0.75pt]    {$1$};
% Text Node
\draw (426.53,57.94) node [anchor=north west][inner sep=0.75pt]    {$0$};
% Text Node
\draw (219.33,175.08) node [anchor=north west][inner sep=0.75pt]    {$0$};
% Text Node
\draw (192.6,145.37) node [anchor=north west][inner sep=0.75pt]    {$0$};
% Text Node
\draw (243.33,145.37) node [anchor=north west][inner sep=0.75pt]    {$0$};
% Text Node
\draw (215.9,118.43) node [anchor=north west][inner sep=0.75pt]    {$0$};
% Text Node
\draw (166.55,118.43) node [anchor=north west][inner sep=0.75pt]    {$0$};
% Text Node
\draw (266.63,119.43) node [anchor=north west][inner sep=0.75pt]    {$0$};
% Text Node
\draw (189.17,88.72) node [anchor=north west][inner sep=0.75pt]    {$2$};
% Text Node
\draw (241.95,88.72) node [anchor=north west][inner sep=0.75pt]    {$1$};
% Text Node
\draw (214.53,55.94) node [anchor=north west][inner sep=0.75pt]    {$3$};

\end{tikzpicture}
\caption{Birth-death function $\BD \phi: \Int{P} \to \Z$ (left) associated to the free presentation from Figure~\ref{fig:presentation} and its persistence diagram $\partial \BD \phi: \Int{P} \to \Z$ (right)}
\label{fig:pres-zb}
\end{center}
\end{figure}
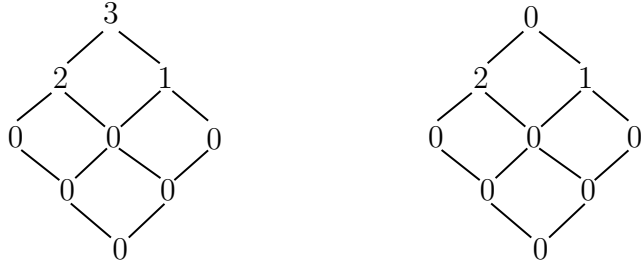

See Figure~\ref{fig:pres-zb} for an example of a birth-death function.

\begin{prop}\label{prop:BD_module}
Let $\phi : F \Rightarrow M$ be a presentation of $M$ and $f : M \rightleftarrows N : g$ be a module-preserving morphism. Then, $\psi := \phi \circ g : F \circ g \Rightarrow N$ is a free presentation of $N$.
Further,  $\BD \psi = (\Int g)^\sharp \BD \phi$.
\end{prop}
\begin{proof}
See \cite[Lemma 4.12]{GulenMcCleary} for a proof of why $\psi : F  \circ g \Rightarrow N$  is a free presentration of~$N$.
Both $N$ and $\psi$ are the restrictions of $M$ and $\phi$ to $\image g$.
This means that the pullback diagrams defining $\BD \phi \big[ g(x), g(y) \big]$ and $\BD \psi [x,y]$ are the same.
\end{proof}

%%%%%%%%%%%%%%%%%%%%%%%%%%%%%%%%%%%%%%%%%%%%%%
\subsection{Persistence Diagram}
%%%%%%%%%%%%%%%%%%%%%%%%%%%%%%%%%%%%%%%%%%%%%%
We interpret $\BD \phi [a,b]$ as generators born by $a$ and dead by $b$.
However, the persistence diagram of $M$ should capture generators born at $a$ and dead at $b$.
We achieve this by taking the combinatorial derivative.
Some care is required so that the persistence diagram of the module~$M$ is independent of the choice of a free presentation~$\phi$.

It is now convenient to say two functions $m, n : \Int P \to \Z$ are \emph{equivalent}, written $m \sim n$, if $m [a,b] = n[a,b]$, for all $a < b$.
That is, $m \sim n$ if they are the same up to the diagonal.

    \begin{prop}\label{prop:dgm_well_defined}
        Let $\phi : F \Rightarrow M$ and $\psi : G \Rightarrow M$ be two free
        presentations of a $P$-module~$M$.
        Then, $\partial \BD \phi \sim \partial \BD \psi$.
    \end{prop}
    \begin{proof}
        See \cite[Corollary 5.8]{GulenMcCleary} for a proof.
    \end{proof}

The above proposition allows us to define persistence diagrams for modules independent of the choice of free presentations.

    \begin{defn}\label{defn:presentation}
    The \define{persistence diagram} of a $P$-module $M$ is the equivalence
    class of the function $\partial \BD \phi$, for any free presentation $\phi : F \Rightarrow M.$
    \end{defn}

See Figure~\ref{fig:pres-zb} for an example of a persistence diagram.

    \begin{prop}[\cite{GulenMcCleary}]
    \label{prop:dgm_galois}
        Let $\phi : F \Rightarrow M$ be a presentation of a $P$-module $M$ 
        and $f : M \rightleftarrows N : g$ a module-preserving morphism.
        Let $\psi := \phi \circ g$ be the induced free presentation of $G$.
        Then, $\partial \BD \psi = (\Int f)_\sharp \big( \partial \BD \phi F \big)$.  
    \end{prop}
    \begin{proof}
    The statement follows from Proposition~\ref{prop:BD_module} and Theorem~\ref{thm:rota}.
    \end{proof}

%%%%%%%%%%%%%%%%%%%%%%%%%%%%%%%%%%%%%%%%%%%%%%
%%%%%%%%%%%%%%%%%%%%%%%%%%%%%%%%%%%%%%%%%%%%%%
\section{Equivalence of Methods}
\label{sec:equivalence}
%%%%%%%%%%%%%%%%%%%%%%%%%%%%%%%%%%%%%%%%%%%%%%
%%%%%%%%%%%%%%%%%%%%%%%%%%%%%%%%%%%%%%%%%%%%%%

We now have two methods for arriving at the persistence diagram of a (co)filtration.
The first makes use of (co)cycles and (co)boundaries.
The second is a purely algebraic construction on its (co)homology persistence module.
We now show the two produce equivalent persistence diagrams.

%%%%%%%%%%%%%%%%%%%%%%%%%%%%%%%%%%%%%%%%%%%%%%
\subsection{Kernel Function}
%%%%%%%%%%%%%%%%%%%%%%%%%%%%%%%%%%%%%%%%%%%%%%
Fix a persistence module $M : P \to \Vect$.
Above, we define the persistence diagram of $M$ as the M\"obius inversion of the birth-death function associated to any presentation of $M$.
To show that is well defined, in other words, for the proof of Proposition \ref{prop:dgm_well_defined}, McCleary and G\"ulen use the M\"obius inversion of an intermediate function  called the kernel function.

\begin{defn}
The \define{kernel function} of $M$ is the function $\ker M: \Int P \to \Z$ defined as $\ker M [a,b] = \dim  \ker M (a \leq b)$.
\end{defn}

\begin{figure}
\begin{center}
% !TEX root = main.tex

\tikzset{every picture/.style={line width=0.75pt}} %set default line width to 0.75pt        

\begin{tikzpicture}[x=0.75pt,y=0.75pt,yscale=-1,xscale=1]
%uncomment if require: \path (0,260); %set diagram left start at 0, and has height of 260

%Straight Lines [id:da1293556665037825] 
\draw    (137,173) -- (170,135) ;
%Straight Lines [id:da7350135388052794] 
\draw    (116,172) -- (81,136) ;
%Straight Lines [id:da9947720037844574] 
\draw    (80,109) -- (112,74) ;
%Straight Lines [id:da19064428364535213] 
\draw    (169,109) -- (139,74) ;
%Straight Lines [id:da9397453688249273] 
\draw    (324,75) -- (306,91) ;
%Straight Lines [id:da12952153281805967] 
\draw    (298,106) -- (280.66,120.04) ;
%Straight Lines [id:da5747406792819012] 
\draw    (353,105) -- (333,123) ;
%Straight Lines [id:da6228674862031855] 
\draw    (326,134) -- (310,149) ;
%Straight Lines [id:da6932310746662698] 
\draw    (355,162) -- (337,179) ;
%Straight Lines [id:da5583241735256144] 
\draw    (377,135) -- (360,150) ;
%Straight Lines [id:da41615788129376674] 
\draw    (352.64,91.02) -- (332,75) ;
%Straight Lines [id:da41794652426687295] 
\draw    (328,179) -- (307.77,161.8) ;
%Straight Lines [id:da383481786753729] 
\draw    (302,150) -- (282.77,134.8) ;
%Straight Lines [id:da6358066397104591] 
\draw    (325,123) -- (304.77,105.8) ;
%Straight Lines [id:da6006793198932561] 
\draw    (377,120) -- (359,105) ;
%Straight Lines [id:da07069541937950463] 
\draw    (353,149) -- (332,134) ;
%Straight Lines [id:da32916043680706286] 
\draw    (536,75) -- (518,91) ;
%Straight Lines [id:da41131331728175646] 
\draw    (510,106) -- (492.66,120.04) ;
%Straight Lines [id:da19122076685705292] 
\draw    (565,105) -- (545,123) ;
%Straight Lines [id:da852272214176335] 
\draw    (538,134) -- (522,149) ;
%Straight Lines [id:da10631214916506071] 
\draw    (567,162) -- (549,179) ;
%Straight Lines [id:da43107398565094845] 
\draw    (589,135) -- (572,150) ;
%Straight Lines [id:da8809421755882967] 
\draw    (564.64,91.02) -- (544,75) ;
%Straight Lines [id:da5636371521203614] 
\draw    (540,179) -- (519.77,161.8) ;
%Straight Lines [id:da04556327797373516] 
\draw    (514,150) -- (494.77,134.8) ;
%Straight Lines [id:da9776367403804185] 
\draw    (537,123) -- (516.77,105.8) ;
%Straight Lines [id:da4224221249165634] 
\draw    (589,120) -- (571,105) ;
%Straight Lines [id:da4617777915316812] 
\draw    (565,149) -- (544,134) ;

% Text Node
\draw (539.33,178.08) node [anchor=north west][inner sep=0.75pt]    {$0$};
% Text Node
\draw (512.6,148.37) node [anchor=north west][inner sep=0.75pt]    {$0$};
% Text Node
\draw (563.33,148.37) node [anchor=north west][inner sep=0.75pt]    {$0$};
% Text Node
\draw (535.9,121.43) node [anchor=north west][inner sep=0.75pt]    {$0$};
% Text Node
\draw (486.55,121.43) node [anchor=north west][inner sep=0.75pt]    {$0$};
% Text Node
\draw (586.63,121.43) node [anchor=north west][inner sep=0.75pt]    {$0$};
% Text Node
\draw (509.17,91.72) node [anchor=north west][inner sep=0.75pt]    {$2$};
% Text Node
\draw (561.95,91.72) node [anchor=north west][inner sep=0.75pt]    {$1$};
% Text Node
\draw (534.53,59.94) node [anchor=north west][inner sep=0.75pt]    {$0$};
% Text Node
\draw (327.33,178.08) node [anchor=north west][inner sep=0.75pt]    {$0$};
% Text Node
\draw (300.6,148.37) node [anchor=north west][inner sep=0.75pt]    {$0$};
% Text Node
\draw (351.33,148.37) node [anchor=north west][inner sep=0.75pt]    {$0$};
% Text Node
\draw (323.9,121.43) node [anchor=north west][inner sep=0.75pt]    {$0$};
% Text Node
\draw (274.55,121.43) node [anchor=north west][inner sep=0.75pt]    {$0$};
% Text Node
\draw (374.63,122.43) node [anchor=north west][inner sep=0.75pt]    {$0$};
% Text Node
\draw (297.17,91.72) node [anchor=north west][inner sep=0.75pt]    {$2$};
% Text Node
\draw (349.95,91.72) node [anchor=north west][inner sep=0.75pt]    {$1$};
% Text Node
\draw (322.53,59.94) node [anchor=north west][inner sep=0.75pt]    {$0$};
% Text Node
\draw (120.09,54.86) node [anchor=north west][inner sep=0.75pt]  [font=\normalsize]  {$0$};
% Text Node
\draw (69.27,112.36) node [anchor=north west][inner sep=0.75pt]  [font=\normalsize]  {$\mathbb{F}^{2}$};
% Text Node
\draw (164.31,113.22) node [anchor=north west][inner sep=0.75pt]  [font=\normalsize]  {$\mathbb{F}$};
% Text Node
\draw (120.67,177.4) node [anchor=north west][inner sep=0.75pt]  [font=\normalsize]  {$0$};

\end{tikzpicture}
\caption{The persistence module $M: P \to \Vect$ (left), its kernel function $\ker M: \Int{P} \to \Z$ (middle), and its persistence diagram $\partial \ker M: \Int{P} \to \Z$ (right). }
\label{fig:kernel}
\end{center}
\end{figure}
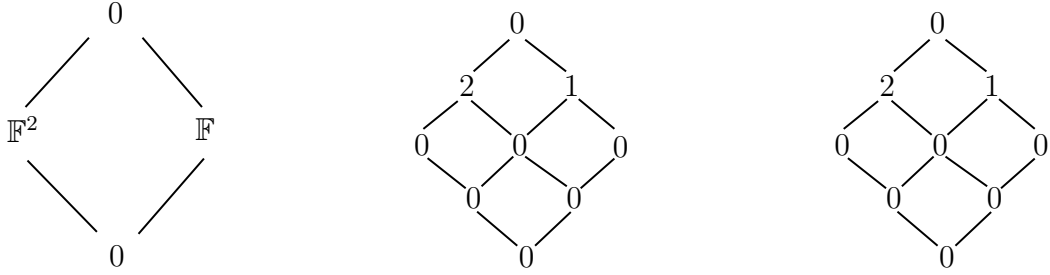

See Figure \ref{fig:kernel} for an example of a kernel function.

\begin{lem}
\label{lem:kernel}
Let~$\BD \phi :  \Int P \to \Z$ be the birth-death function of any presentation~$\phi : F \Rightarrow M$.
Then, $\partial \BD \phi \sim  \partial \ker M$.
\end{lem}
\begin{proof}
See proof of \cite[Proposition~5.7]{GulenMcCleary}.
\end{proof}

%%%%%%%%%%%%%%%%%%%%%%%%%%%%%%%%%%%%%%%%%%%%%%
\subsection{Filtrations}
%%%%%%%%%%%%%%%%%%%%%%%%%%%%%%%%%%%%%%%%%%%%%%

We have two persistence diagrams associated to a filtration $F : P \to \Delta K$.
The first $\partial \BD_\ast F$ is via cycles and boundaries  (Definition~\ref{defn:cycles_boundaries}).
The second $\partial \BD  \phi$ is via its persistent homology module $H_\ast | F |: P \to \Vect$ obtained by applying singular homology to the underlying spaces~$|F(a)|$ (Definition~\ref{defn:presentation}).

\begin{prop}[\cite{mccleary_communication}]
\label{prop:filtration_equivalence}
Let $F : P \to \Delta K$ be a filtration and
consider its two persistence diagrams $\partial \BD_\ast F$ and $\partial \BD \phi$, where $\phi$ is any free presentation of the homology persistence module $H_\ast |F|$.
Then, $\partial \BD_\ast F \sim \partial \BD \phi$.
\end{prop}
\begin{proof}
By Lemma~\ref{lem:kernel}, $\partial \BD \phi \sim \partial \ker H_\ast |F|$.
We show $\partial \ker H_\ast |F| \sim  \partial \BD_\ast F$.

Let $B_\ast F : \Int P \to \Z$ be the function $B_\ast F  [a,b] = \dim B_\ast F(a)$.
We have
    $$\ker H_\ast |F| [a,b] = 
    \dim \left( \frac{Z_\ast F(a) \cap B_\ast F(b)}{B_\ast F(a)} \right) = \BD_\ast F [a,b] - B_\ast F [a,b].$$
By linearity of the M\"obius inversion, we have
    $$\partial \ker H_\ast |F| = \partial \BD_\ast F - \partial B_\ast F.$$
To finish the proof, we show $\partial B_\ast F \sim 0$.
Let $f : P \leftrightarrows \Int P : g$ be the Galois connection $f(a) = [a,a]$ and $g[a,b] = a$.
Consider the function $\beta : P \to \Z$ defined as $\beta(a) = \dim B_\ast F(a)$.
Note $B_\ast F = g^\sharp \beta$.
By Theorem~\ref{thm:rota}, $\partial B_\ast F = f_\sharp \partial \beta$.
Since $f^{-1}[a,b] = \emptyset$ for every $a < b$,
we  have $\partial B_\ast F [a,b] = 0$ as desired.
\end{proof}

%%%%%%%%%%%%%%%%%%%%%%%%%%%%%%%%%%%%%%%%%%%%%%
\subsection{Cofiltrations}
%%%%%%%%%%%%%%%%%%%%%%%%%%%%%%%%%%%%%%%%%%%%%%
We have two persistence diagrams associated to a cofiltration $F : P \to \nabla K$.
The first~$\partial \BD^\ast F$ is via cocycles and coboundaries (Definition~\ref{defn:cocyles_coboundaries}).
The second $\partial \BD \phi$ is via its persistent cohomology module $H^\ast_\com | F |: P \to \Vect$ obtained by applying singular cohomology with compact support to the underlying spaces $|F(a)|$ (Definition~\ref{defn:presentation}).

\begin{prop}
\label{prop:cofiltration_equivalence}
Let $F : P \to \nabla K$ be a cofiltration and
consider its two persistence diagrams $\partial \BD^\ast F$ and $\partial \BD \phi$, where $\phi$ is any free presentation of the cohomology persistence module $H^\ast_\com |F|$.
Then, $\partial \BD^\ast F \sim \partial \BD \phi$.
\end{prop}
\begin{proof}
The proof is the same as that of Proposition \ref{prop:filtration_equivalence}.
Simply replace cycles and boundaries with cocyles and coboundaries.
\end{proof}

\begin{ex}
In Figures~\ref{fig:presentation}, \ref{fig:pres-zb}, and \ref{fig:kernel}, we chose the persistence module $M$ to be given by $H^1_\com |F|$, the persistent cohomology module of the cofiltration $F$ from Figure~\ref{fig:mob-ph-ex}. Indeed, the persistent cohomology diagrams in Figures~\ref{fig:pres-zb} and \ref{fig:kernel} are identical, and the persistent cohomology diagram in Figure~\ref{fig:mob-ph-ex} differs from the other two only along the diagonal.
\end{ex}

%%%%%%%%%%%%%%%%%%%%%%%%%%%%%%%%%%%%%%%%%%%%%%
%%%%%%%%%%%%%%%%%%%%%%%%%%%%%%%%%%%%%%%%%%%%%%
\section{Poincar\'e Duality}
\label{sec:duality}
%%%%%%%%%%%%%%%%%%%%%%%%%%%%%%%%%%%%%%%%%%%%%%
%%%%%%%%%%%%%%%%%%%%%%%%%%%%%%%%%%%%%%%%%%%%%%

We now assemble the three definitions of a persistence diagram into a Poincar\'e duality between filtrations and cofiltrations of a manifold.

\begin{figure}
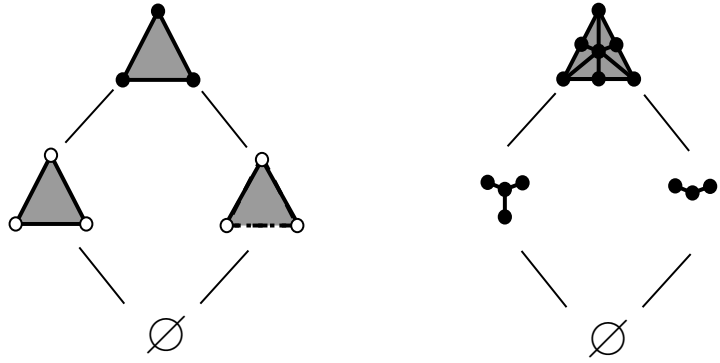

\begin{center}
\include{figures-tex/dual-filt.tex}
\caption{A cofiltration $F: P \to \nabla K$ (left) and its dual filtration $G: P \to \Delta L$ (right)}
\label{fig:dual-filt}
\end{center}
\end{figure}

Every cofiltration $F : P \to \nabla K$ has a \emph{dual filtration} $G : P \to \Delta L$ as follows.
Let~$L$ be the barycentric subdivision of $K$.
Recall, simplices of $L$ are chains $\sigma_0 < \cdots < \sigma_n$ of simplices in $K$ and the face relation is the subchain relation.
For $a \in P$, let $G(a) \subseteq L$ be all chains 
$\sigma < \cdots $ starting with a simplex $\sigma$ in $F(a)$. 
Since $F(a) \subseteq K$ is a supcomplex, every subchain of such a chain also starts in $F(a)$.
This makes $G(a)$ a subcomplex.
The open star of $G(a)$ in $L$, denoted $\st G(a)$,
is the set of all simplices in $L$ that has a face in $G(a)$.
In other words, $\st G(a)$ are all chains $\sigma < \cdots$ in $K$
with $\sigma \in F(a)$.
This means $\st G(a)$ is the barycentric subdivision of $F(a)$.
An important property of the open star of any subcomplex $A \subseteq L$
is that $|A|$ is a deformation retract of $| \st A |$~\cite[Lemma 70.1]{munkres}.
With this, we conclude that $|G(a)|$ is a deformation retract of~$|F(a)|$. 
See Figure~\ref{fig:dual-filt} for an example of a dual filtration.

\begin{thm}\label{thm:duality_one}
Let $K$ be a triangulation of a compact, $\field$-orientable, $m$-manifold.
For any cofiltration $F : P \to \nabla K$, let $G : P \to \Delta L$ be its dual filtration.
Then $\partial \BD^i F \sim \partial \BD_{m-i} G$.
\end{thm}
\begin{proof}
Let $H^i_\com |F|$ be the persistence module obtained by applying singular cohomology with compact support to $F$.
Let $H_{m-i} |F|$ be the persistence module obtained by applying singular homology to $F$.
By Theorem~\ref{thm:poincare}, $H^i_\com |F| \cong H_{m-i} |F|$.

Let $H_{m-i} |G|$ be the persistence module obtained by applying singular homology to $G$.
Since  $|F(a)|$ is homotopy equivalent to $|G(a)|$, for every $a \in P$, $H_{m-i} |G| \cong H_{m-i} |F|$.

These isomorphisms combined with Propositions~\ref{prop:filtration_equivalence} and~\ref{prop:cofiltration_equivalence} imply the desired result.
\end{proof}

\begin{figure}
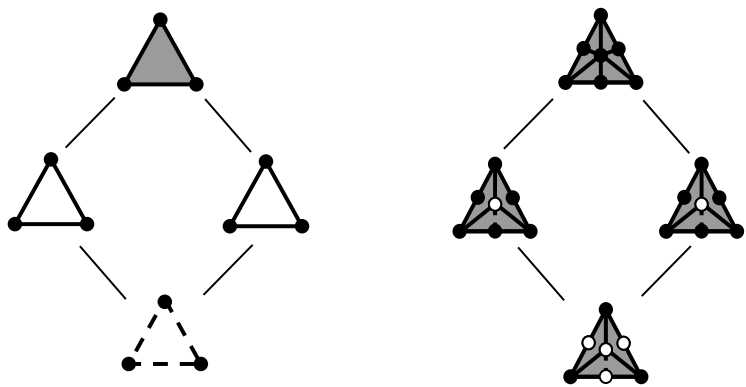

\begin{center}
\include{figures-tex/dual-cofilt.tex}
\caption{A filtration $F: P \to \Delta K$ (left) and its dual cofiltration $G: P \to \nabla L$ (right)}
\label{fig:dual-cofilt}
\end{center}
\end{figure}

Similarly, every filtration $F : P \to \Delta K$ has a \emph{dual cofiltration} $G : P \to \nabla L$ as follows.
Let~$L$ be the barycentric subdivision of $K$.
For $a \in P$, let $G(a) \subseteq L$ be all chains 
$\sigma < \cdots$ starting with a  simplex $\sigma$ in  $F(a)$.
Since $F(a) \subseteq K$ is a subcomplex, every supchain of such a chain starts with a simplex in $F(a)$.
This makes $G(a)$ a supcomplex.
Note that the open star of the subdivision of $F(a)$ in $L$ is $G(a)$.
This means that $|F(a)|$ is a deformation retract of $|G(a)|$.
See Figure~\ref{fig:dual-cofilt} for an example of a dual cofiltration.

\begin{thm}\label{thm:duality_two}
Let $K$ be a triangulation of a compact, $\field$-orientable, $m$-manifold.
For any filtration $F : P \to \Delta K$, let $G : P \to \nabla L$ be its dual cofiltration.
Then $\partial \BD_{m-i} F \sim \partial \BD^{i} G$.
\end{thm}
\begin{proof}
The proof is the proof of Theorem~\ref{thm:duality_one} but with the roles
of $F$ and $G$ reversed. 
\end{proof}

\begin{comment}

\begin{figure}
\begin{center}
\include{figures-tex/torus_filt_cofilt}
\caption{A cofiltration $F$ of the triangulated torus, $T$, along with its dual filtration $G$. \todo{Label/color to show identifications}}
\label{fig:torus-filt-cofilt}
\end{center}
\end{figure}

\begin{figure}
\begin{center}
\include{figures-tex/torus_pds}
\caption{Persistence diagrams for the cofiltration $F$ and its dual filtration $G$ from Figure \ref{fig:torus-filt-cofilt}.}
\label{fig:torus-pds}
\end{center}
\end{figure}

\begin{ex}
See Figures \ref{fig:torus-filt-cofilt} and \ref{fig:torus-pds} for an illustration of Theorem \ref{thm:duality_one}. In Figure \ref{fig:torus-filt-cofilt}, we give a cofiltration $F$ of the triangulated torus $T$ and its dual filtration $G$ of the barycentric subdivison $\mathcal{B}(T)$ of $T$. In Figure \ref{fig:torus-pds}, we compute persistence diagrams for $F$ and $G$: indeed, $\partial \BD^i F \sim \partial \BD_{2-i}G$ for $i = 0, 1, 2$.
\end{ex}
\end{comment}

%%%%%%%%%%%%%%%%%%%%%%%%%%%%%%%%%%%%%%%%%%%%%%
%%%%%%%%%%%%%%%%%%%%%%%%%%%%%%%%%%%%%%%%%%%%%%
\section{Statements and Declarations}
%%%%%%%%%%%%%%%%%%%%%%%%%%%%%%%%%%%%%%%%%%%%%%
%%%%%%%%%%%%%%%%%%%%%%%%%%%%%%%%%%%%%%%%%%%%%%
\paragraph{Conflict of Interest}
The authors have no financial or proprietary interests in any material discussed in this article.

%%%%%%%%%%%%%%%%%%%%%%%%%%%%%%%%%%%%%%%%%%%%%%
\bibliographystyle{plain}
\bibliography{references}
%%%%%%%%%%%%%%%%%%%%%%%%%%%%%%%%%%%%%%%%%%%%%%

\end{document}